\documentclass[12pt]{article}
\usepackage{latexsym}
\usepackage{amssymb}
\usepackage{amsfonts}
\usepackage{amsmath}
\usepackage{enumerate}
\begin{document}
\title{Eternal Domination in Trees}
\author{William F. Klostermeyer$^{\dagger}$ and Gary
MacGillivray$^{\ddagger}$
\vspace*{0.07in}\\
      $^{\dagger}$ School of Computing\\
      University of North Florida \\ Jacksonville, FL 32224-2669
\vspace*{0.07in}\\
     $^{\ddagger}$ Dept. of Mathematics and Statistics\\
       University of Victoria\\ Victoria, Canada\\}

\date{}
\maketitle

\begin{abstract}
Mobile guards on the vertices of a graph are used to defend the
graph against an infinite sequence of attacks on vertices. A guard
must move from a neighboring vertex to an attacked vertex (we assume
attacks happen only at vertices containing no guard). More than one
guard is allowed to move in response to an attack. The $m$-eternal
domination number is the minimum number of guards needed to defend
the graph. We characterize the trees achieving several upper and
lower bounds on the $m$-eternal domination number.
\end{abstract}
\thispagestyle{empty}

\newtheorem{theorem}{Theorem}[section]
\newtheorem{other}{Theorem}
\newtheorem{lemma}[theorem]{Lemma}
\newtheorem{prop}[theorem]{Proposition}
\newtheorem{cor}[theorem]{Corollary}
\newtheorem{conj}{Conjecture}
\newtheorem{prob}[conj]{Problem}
\newtheorem{question}[conj]{Question}
\newtheorem{fact}[theorem]{Fact}
\newtheorem{proposition}[theorem]{Proposition}
\newtheorem{corollary}[theorem]{Corollary}

\noindent\textbf{Keywords:\hspace{0.1in}} dominating set, eternal
dominating set, connected dominating set, independent set,
neo-colonization.


\section{Introduction}

Let $G=(V,E)$ be a graph with $n$ vertices. Several recent papers
have considered problems associated with using mobile guards to
defend $G$ against an infinite sequence of attacks; see for instance
\cite{ABB, BCG2, GHH, GK, KM, KM2, KM6, KM7}.

Denote the open and closed neighborhoods of a vertex $x\in V$ by
$N(x)$ and $N[x]$, respectively. That is, $N(x)=\{v|xv\in E\}$ and
$N[x]=N(x)\cup\{x\}$. Further, for $S\subseteq V$, let
$N(S)=\bigcup_{x\in S}N(x)$. For any $X\subseteq V$ and $x\in X$, we
say that $v\in V-X$ is an \emph{external private neighbor }of $x$
\emph{with respect to} $X$ if $v$ is adjacent to $x$ but to no other
vertex in $X$; we sometimes simply say that $v$ is a
private neighbor of $x$. The set of all such vertices $v$ is the
\emph{external private neighborhood of }$x$ \emph{with respect to
}$X$.

A \emph{dominating set }of graph $G$ is a set $D\subseteq V$ with
the property that for each $u\in V-D$, there exists $x\in D\
$adjacent to $u$. A dominating set $D$ is a \emph{connected
dominating set} if the subgraph $G[D]$ induced by $D$ is connected.
The minimum cardinality amongst all dominating sets of $G$ is the
\emph{domination number} $\gamma(G)$, while the minimum cardinality
amongst all connected dominating sets is the \emph{connected}
\emph{domination number} $\gamma_{c}(G)$. Further background on
domination can be found in \cite{HHS}.

An \emph{independent set} of vertices in $G$ is a set $I\subseteq
V$ with the property that no two vertices in $I$ are adjacent. The
maximum cardinality amongst all independent sets is the
\emph{independence number}, which we denote as $\beta(G)$.

Let $D_{i}\subseteq V,1\leq i$, be a set of vertices with one
guard located on each vertex of $D_{i}$. In this paper, we shall
allow at most one guard to be located on a vertex at any time. The
problems considered in this paper can be modeled as a two-player
game between a \emph{defender} and an \emph{attacker}: the
defender chooses $D_{1}$ as well as each $D_{i},i>1$, while the
attacker chooses the locations of the attacks
$r_{1},r_{2},\ldots$. Note that the location of an attack can be
chosen by the attacker depending on the location of the guards.
Each attack is handled by the defender by choosing the next
$D_{i}$ subject to some constraints that depend on the particular
game. The defender wins the game if they can successfully defend
any series of attacks, subject to the constraints of the game;
the attacker wins otherwise.

In the \emph{eternal dominating set problem}, each $D_{i},i\geq1$,
is required to be a dominating set, $r_{i}\in V$ (assume without
loss of generality $r_{i}\notin D_{i}$), and $D_{i+1}$ is obtained
from $D_{i}$ by moving one guard to $r_{i}$ from a vertex $v\in
D_{i},v\in N(r_{i})$. The smallest size of an eternal dominating
set for $G$ is denoted $\gamma^{\infty}(G)$. This problem was
first studied in \cite{BCG2}.

In the \emph{$m$-eternal dominating set problem}, each
$D_{i},i\geq1$, is required to be a dominating set, $r_{i}\in V$
(assume without loss of generality $r_{i}\notin D_{i}$), and
$D_{i+1}$ is obtained from $D_{i}$ by moving guards to neighboring
vertices. That is, each guard in $D_{i}$ may move to an adjacent
vertex. It is required that $r_{i}\in D_{i+1}$. The smallest size of an
$\mathrm{m}$-eternal dominating set for $G$ is denoted
$\gamma _{\mathrm{m}}^{\infty}(G)$. This \textquotedblleft
all-guards move\textquotedblright\ version of the problem was
introduced in \cite{GHH}. It is clear that
$\gamma^{\infty}(G)\geq\gamma_{\mathrm{m}}^{\infty}(G)\geq\gamma(G)$ for all graphs $G$.

We say that a vertex is \emph{protected} if there is a guard on
the vertex or on an adjacent vertex. We say that an attack at $v$
is \emph{defended} if we send a guard to $v$.

Our objective in this paper is to describe the trees that achieve
some upper and lower bounds on the $m$-eternal domination number.
In Sections 3 through 6 we describe the trees for which equality holds in each of the following:
\begin{enumerate}
\item[(a)] $\gamma_{\mathrm{m}}^{\infty}(T) \leq \gamma_c(T)+1$;
\item[(b)] $\gamma(T)  \leq \gamma_{\mathrm{m}}^{\infty}(T)$;
\item[(c)] $\gamma_{\mathrm{m}}^{\infty}(T) \leq 2\gamma(T)$; and
\item[(d)] $\gamma_{\mathrm{m}}^{\infty}(T) \leq \beta(T)$.
\end{enumerate}

\section{Terminology and Background\label{SecTrees}}

A \emph{neo-colonization} is a partition
$\{V_{1},V_{2},\ldots,V_{t}\}$ of the vertex set of graph $G$ such
that each $G[V_{i}]$ is a connected graph. A part $V_{i}$ is
assigned a weight $\omega(V_{i})=1$ if it induces a clique and
$\omega(V_{i})=1+\gamma_{c}(G[V_{i}])$ otherwise. Then
$\theta_{c}(G)$ is the minimum total weight of any
neo-colonization of $G$, and is called the \emph{clique-connected
cover number} of $G$. Goddard et al. \cite{GHH} defined this
parameter and proved that
$\gamma_{\mathrm{m}}^{\infty}(G)\leq\theta _{c}(G)$.
For $X \subseteq V$ we will write $\theta_c(X)$ as shorthand for
$\theta_c(G[X])$.
Klostermeyer
and MacGillivray then proved the next result, which is key to what
follows in this paper.

\begin{theorem}
\emph{\cite{KM2}}\label{neo-c} For any tree $T$, $\theta_{c}(T)=\gamma
_{\mathrm{m}}^{\infty}(T)$.
\end{theorem}

It follows from the previous theorem that
$\gamma_{\mathrm{m}}^{\infty}(G)\leq\gamma _{c}(G)+1$.

\medskip

The following property of neo-colonizations will be useful.

\begin{proposition}\label{nosingles}
Let $T$ be a tree with at least two vertices.  Then there is a
neo-colonization of minimum weight in which every part has size at
least two.
\end{proposition}

\noindent{\it Proof.} Suppose that  $\Pi = \{V_1, V_2, \ldots, V_k,
\{x\}\}$ is a minimum weight neo-colonization of $T$.  Without loss
of generality, $\{x\}$ is adjacent to a vertex of $V_k$.  Since
the subgraph of $T$ induced by $V_k$ is connected, $x$ is adjacent
to exactly one vertex of $V_k$.  The neo-colonization $\Pi^\prime
= \{V_1, V_2, \ldots, V_k \cup \{x\}\}$ has minimum weight and fewer
parts of size one.  Applying this argument repeatedly, one arrives
at the desired neo-colonization. $\Box$

\medskip

A \emph{stem }of a tree $T$ is a vertex of degree at least two
that is adjacent to a leaf.
A vertex of $T$ that
is not a leaf is called an \emph{internal} vertex. A tree is a
{\it star} if it is isomorphic to $K_{1, m}, m \geq 1$.

We partition the internal vertices of $T$ into \emph{loners},
\emph{weak stems }and \emph{strong stems }depending on whether
they are adjacent to no, exactly one or at least two leaves.
Denote the set of leaves of $T$ by $L(T)$ and let $\ell=|L(T)|$.
Obviously, $\gamma_{c}(T)=n-\ell$, the number of \emph{internal vertices}, for any tree $T$ of order
$n\geq3$.

The \emph{eccentricity} of a vertex in a graph is its maximum distance from any other vertex.  A vertex of maximum eccentricity in a tree is a leaf which is an end vertex of a longest path.  Leaves of maximum eccentricity and the stems to which they are adjacent play an important role in some of our proofs. We use $\mathit{deg}(v)$ to denote the degree of vertex $v$.





%

\section{Trees with  $\gamma_{\mathrm{m}}^\infty = \gamma_c + 1$}

It turns out to be  easier to describe the trees for which  $\gamma_{\mathrm{m}}^\infty < \gamma_c + 1$. The results below make it possible to look at a tree and determine if $\gamma_{\mathrm{m}}^\infty < 1 + \gamma_c$.  They do not, however, give much structural information on the trees for which the inequality holds.  Finding such results is an open problem.

\begin{proposition}
 A tree $T$ has $\gamma_{\mathrm{m}}^\infty < 1 + \gamma_c$ if and only if $T$ has a spanning forest consisting of $r\ K_2$'s and trees $T_1, T_2, \ldots, T_k$ on at least three vertices
such that at least $k$ loners of $T$ are leaves of the $k+r$ trees in the collection.
\label{kLeaves}
\end{proposition}

\noindent{\it Proof.} Suppose $\gamma_{\mathrm{m}}^\infty < 1 + \gamma_c$.  By Lemma \ref{nosingles}, there is a minimum weight neo-colonization $\Pi = \{V_1, V_2, \ldots, V_{k+r}\}$ in which there are no parts of size one.  Without loss of generality, $V_1, V_2, \ldots, V_k$ each have size at least three, and the remaining parts each have size two.  Let $T_i = T[V_i],\ 1 \leq i \leq k+r$.  Every internal vertex of some $T_i$ is an internal vertex of  $T$.  Since there are no parts of size one, every leaf of $T$ belongs to the same part as the stem to which it is adjacent. In particular, no stem of $T$ is a leaf of any tree $T_i$ that is not isomorphic to $K_2$.  Each of the trees $T_j,\ k+1 \leq j \leq k+r$ contains at least one internal vertex of $T$.  The weight of $\Pi$ is $r + k$ plus the total number of internal vertices of the trees $T_i,\ 1 \leq i \leq k$.  All but $k$ units of this quantity are accounted for by internal vertices of $T$.  Since $\gamma_c$ equals the number of internal vertices of $T$, it follows that at least $k$ loners of $T$ appear as leaves of the trees in the collection.

On the other hand, suppose that $T$ has a spanning forest consisting of trees $T_1, T_2, \ldots, T_k$ on at least three vertices and $r$ trees $T_{k+1}, T_{k+2} \ldots, T_{k+r}$, each isomorphic to $K_2$, such that at least $k$ loners of $T$ are leaves of the $k+r$ trees in the collection.  Let $V_i = V[T_i],\ 1 \leq i \leq k+r$ and consider the neo-colonization $\Pi = \{V_1, V_2, \ldots, V_{k+r}\}$.  The weight of $\Pi$ is $r + k$ plus the total number of internal vertices of the trees $T_i,\ 1 \leq i \leq k$.  Since each of the trees which are isomorphic to $K_2$ contains at least one internal vertex of $T$, and every leaf of $T$ belongs to the same part as the stem to which it is adjacent, this quantity is at most $\gamma_c$, the number of internal vertices of $T$.
$\Box$

\begin{corollary}
Let $T$ be a tree.  Then $\gamma_{\mathrm{m}}^\infty < 1 + \gamma_c$ if and only if
there exists a set of edges whose deletion creates a spanning forest
consisting of $r\ K_2$'s and trees $T_1, T_2, \ldots, T_k$ on at least three vertices
such that at least $k$ loners of $T$ are leaves of the $k+r$ trees in the collection.
\end{corollary}

\noindent{\it Proof.} The implication that if the condition holds then $\gamma_{\mathrm{m}}^\infty < 1 + \gamma_c$ follows from Proposition \ref{kLeaves}.

Suppose $\gamma_{\mathrm{m}}^\infty < 1 + \gamma_c$.
By Proposition \ref{kLeaves} the tree $T$ has a spanning forest consisting of $r\ K_2$'s and trees $T_1, T_2, \ldots, T_k$ on at least three vertices
such that at least $k$ loners of $T$ are leaves of the $k+r$ trees in the collection.  The set of edges to delete consists of the edges of $T$ with ends in different subtrees.
$\Box$

\medskip
We note, in particular, that if every internal vertex is a weak stem then the corollary holds with $k=0$.

\section{Trees with $\gamma_{\mathrm{m}}^\infty = \gamma$}

Informally, the corona of a graph $G$ is the graph obtained by joining a new vertex of degree one to each vertex of $G$.  Formally, if $G$ has vertex set $V(G) = \{v_1, v_2, \ldots, v_n\}$, then $\mathit{corona}(G)$ is the graph with the $2n$ vertices
$V(\mathit{corona}(G)) = \{v_1, v_2, \ldots, v_n\} \cup \{v_1^\prime, v_2^\prime, \ldots, v_n^\prime\}$ and edges
$E(\mathit{corona}(G)) = E(G) \cup \{v_i v_i^\prime: 1 \leq i \leq n\}$.  A graph $H$ \emph{is a corona} if it is the corona of some graph $G$.  It is known that a connected graph has domination number $\frac{|V|}{2}$ if and only if it is either $C_4$ or a corona \cite{Fink,Payan}.  The trees with $\gamma_{\mathrm{m}}^\infty = \gamma$ turn out to be exactly the coronas of trees.

\begin{lemma}\label{EqualsGamma}
Let $T$ be a tree for which $\gamma_{\mathrm{m}}^{\infty }(T) = \gamma(T)$.  For each minimum weight neo-colonization $\Pi = \{V_1, V_2, \ldots,$ $V_k\}$  of $T$ and $i = 1, 2, \ldots, k$, we have
 $\gamma_{\mathrm{m}}^{\infty }(T[V_i]) = \gamma(T[V_i])$.
\end{lemma}

\noindent{\it Proof.}
Each induced subgraph $T[V_i]$ is connected.
Since $\gamma_{\mathrm{m}}^{\infty }(T)$ is a minimum taken over all neo-colonizations, it must be that for $i = 1, 2, \ldots, k$, the partition $\{V_i\}$ is a neo-colonization of $T[V_i]$ of weight $\gamma_{\mathrm{m}}^{\infty }(T[V_i])$.
Hence,
$$\gamma_{\mathrm{m}}^{\infty }(T)
= \sum_{i=1}^k \omega(V_i)
= \sum_{i=1}^k \gamma_{\mathrm{m}}^{\infty }(T[V_i])
\geq \sum_{i=1}^k \gamma_c(T[V_i])
\geq \sum_{i=1}^k \gamma(T[V_i])
\geq \gamma(T)$$
and the result follows.
 $\Box$

\begin{theorem}\label{DomLowerBound}
A tree $T$ with $n \geq 2$ vertices satisfies $\gamma_{\mathrm{m}}^\infty = \gamma$ if and only if  $T$ is a corona.
\end{theorem}

\noindent{\it Proof.}
Suppose $T$ is a corona.  Then $\gamma(T) = \frac{n}{2}$. It is clear that  there is a neo-colonization of weight $n/2$ in which each part consists of a leaf and its unique neighbor.  Hence $\gamma_{\mathrm{m}}^\infty = \gamma$.

The proof of the converse is by induction on $n$.
The statement is clearly true if $n = 2, 3,4$, the trees $K_2$ and $P_4$ being the only ones with $\gamma_\mathrm{m}^\infty = \gamma$.  Suppose the statement holds for all trees on at least two, and at most $n-1$ vertices, for some $n \geq 5$.  Let $T$ be a tree on $n$ vertices for which $\gamma_\mathrm{m}^\infty = \gamma$.  Then $T$ is not a star.

Let $\Pi = \{V_1, V_2, \ldots, V_k\}$ be a minimum weight neo-colonization of $T$.  By Lemma \ref{EqualsGamma} we have $\gamma_{\mathrm{m}}^{\infty }(T[V_i]) = \gamma(T[V_i])$ for $i = 1, 2, \ldots, k$.  By the induction hypothesis, each tree $T[V_i]$ is a corona.  We claim that, in fact, each is isomorphic to $K_2$.  Otherwise, without loss of generality $T[V_1]$ has at least three vertices and $\omega(V_1) = 1 + \gamma_c(V_1) > \gamma(V_1)$, a contradiction.  This proves the claim.  It follows that $\gamma_{\mathrm{m}}^\infty(T) = k$ and $n = 2k \geq 6$.

Let $u$ be a leaf of maximum eccentricity in $T$, and $v$ be the stem to which it is adjacent.  By the above argument, without loss of generality $V_k = \{u, v\}$.  By Proposition \ref{nosingles}, the vertex $v$ is not adjacent to another leaf besides $u$.  The choice of $u$ now guarantees that $v$ has degree two in $T$.  Hence $T - u - v$ is a tree.  Further,  $\{V_1, V_2, \ldots, V_{k-1}\}$ must be a minimum weight neo-colonization of $T-u-v$.  Thus,
$$\gamma_{\mathrm{m}}^\infty(T-u-v) = k-1 \leq \gamma(T-u-v)  \leq \gamma_{\mathrm{m}}^\infty(T-u-v).$$
It follows that $\gamma(T-u-v) = k-1 = \gamma_{\mathrm{m}}^\infty(T-u-v)$.  By the induction hypothesis, $T-u-v$ is a corona.

Let $w$ be the vertex of $T-u-v$ to which $v$ is adjacent in $T$.  The proof will be complete if we can show that $w$ is a stem.  Otherwise, $w$ is a leaf of the corona $T-u-v$, which has at least four vertices.  Let $s$ be its unique neighbor.  Hence $T$ contains the induced path on four vertices, $u, v, w, s$.  Since the subgraph of $T$ induced by each set $V_i$ is connected and of size two, without loss of generality $V_{k-1} = \{w, s\}$.   Thus,  $T^\prime = T - \{u, v, w\}$ is a tree.  The set $X$ of all internal vertices except $s$ form a dominating set of size $\frac{|V(T^\prime)|-1}{2} = \frac{n-3}{2}$.  But then $X \cup \{v\}$ is a dominating set of $T$ and
$$\gamma(T) \leq \frac{n-1}{2} < \frac{n}{2} = k = \gamma_{\mathrm{m}}^\infty(T),$$
a contradiction.  This completes the proof.  $\Box$

\medskip

An implication of Theorem \ref{DomLowerBound} is that, for trees,
$\gamma_{\rm m}^\infty = \gamma$ only if the clique covering
number equals $\gamma$. This implication holds for all graphs \cite{GHH}. The
converse is not true, even if we restrict the cliques to have size
at least two. For example, $\gamma_{\rm m}^\infty(C_6) =
\gamma(C_6) = 2$, but the clique covering number of $C_6$ equals
three.  The proof of the theorem shows that  trees with
$\gamma_{\rm m}^\infty = \gamma$ have a unique neo-colonization of minimum weight: each part consists of a leaf and the stem to which it is adjacent.


\section{Trees with $\gamma_{\mathrm{m}}^\infty = 2\gamma$}

The following is from \cite{KM6} and we include the proof for
completeness.

\begin{proposition} \cite{KM6} \label{dom-bound}
\label{twice} For any connected graph $G$,
$\gamma_{\mathrm{m}}^{\infty }(G)\leq2\gamma(G)$, and the bound is
sharp for all values of $\gamma(G)$.
\end{proposition}

\noindent\textit{Proof.\hspace{0.1in}} Let $D$ be a minimum
dominating set. Place a guard at each vertex of $D$. For each
vertex $v\in D$, if $v$ has at least one external private
neighbor, pick one of them, say $u$, and place a guard at $u$. It
is easy to see this configuration is an $\mathrm{m}$-eternal
dominating set.

To see that the bound is sharp for $\gamma=1$, consider any star
with at least three vertices. For $\gamma=2$, consider $C_{6}$ and
let $u$ and $v$ be two vertices at distance three apart. Add two
new internally disjoint $u-v$ paths of length three to form the
graph $G$. Obviously, $\{u,v\}$ is a dominating set of $G$. Let $D$
be any minimum dominating set of $G$ with $|D|=3$. Suppose $u\notin D$.
Since $N(u)$ is independent with $|N(u)|=4$, and no two vertices
in $N(u)$ have a common neighbor other than $u$, $D$ does not
dominate $N(u)$, a contradiction. Thus $u\in D$ and similarly
$v\in D$. Without loss of generality say $D=\{u,v,w\}$, where
$w\in N(u)$. Then $D$ cannot repel an attack at a vertex in $N(v) - N(w)$.
It follows that $\gamma_{\mathrm{m}}^{\infty }(G)=4=2\gamma(G)$.

For $\gamma=k\geq3$, consider $C_{3k}$ and let
$\{u_{1},...,u_{k}\}$ be any
minimum dominating set of $C_{3k}$. Note that for each $i$, $d(u_{i}%
,u_{i+1\ (\operatorname{mod}\ 3k)})=3$. For each $i=1,...,k$, add
a new $u_{i}-u_{i+1\ (\operatorname{mod}\ 3k)}$ path of length
three to form $G$. Then $\gamma(G)=k$, but it can be shown similar
to the previous case that no set of $2k-1$ vertices eternally
protects the vertices of $G$. $\Box$

\medskip


We shall now give several characterizations of the trees achieving
the bound in Theorem \ref{dom-bound}.

\begin{proposition}\label{PropIndep}
Let $T$ be a tree.  If $\gamma_{\rm m}^\infty(T) = 2\gamma(T)$
then every minimum dominating set of $T$ is an independent set.
\end{proposition}

\noindent{\it Proof.} We prove the contrapositive.  Suppose there
is a minimum dominating set $D$ and vertices $w, x \in D$ such
that $wx \in E$. Let $\Pi = \{V_1, V_2, \ldots,$ $V_{\gamma}\}$ be a
neo-colonization of $T$ in which, for $i = 1, 2, \ldots, \gamma$,
the subgraph of $T$ induced by $V_i$ is a star containing exactly
one vertex of $D$.  If there exists $i$ such that $|V_i| < 3$,
then $\gamma_{\rm m}^\infty(T) < 2\gamma(T)$. Hence assume that
$|V_i| \geq 3$ for $i = 1, 2, \ldots, \gamma$.  Without loss of
generality, $w \in V_{\gamma-1}$ and $x \in V_{\gamma}$.  Then,
the neo-colonization $\Pi^\prime = \{V_1, V_2, \ldots,$
$V_{\gamma-2}, V_{\gamma-1} \cup V_{\gamma}\}$ has weight at most
$2\gamma-1$ because the connected domination number of
$V_{\gamma-1} \cup V_{\gamma}$ equals two, and hence it
contributes $3 < 4$ to the weight of $\Pi^\prime$. Therefore
$\gamma_{\rm m}^\infty(T) < 2\gamma(T)$. $\Box$

\begin{proposition}\label{Prop2EPN}
Let $T$ be a tree such that $\gamma_{\rm m}^\infty(T) =
2\gamma(T)$.  If $D$ is a minimum dominating set of $T$, then
every $x \in D$ has at least two external private neighbors.
\end{proposition}

\noindent{\it Proof.} We prove the contrapositive.  Suppose first
$D$ is a minimum dominating set of $T$ such that there exists $w
\in D$ with no external private neighbor.  Then $T$ admits a
neo-colonization $\Pi = \{V_1, V_2, \ldots, V_{\gamma-1},$
$\{w\}\}\}$ in which the subgraph induced by $V_i$ is a star
centered at a vertex in $D-\{w\}$.  The weight of $\Pi$ is at most
$2(\gamma-1)+1 < 2 \gamma$.  The argument is similar if $w$ has
exactly one private neighbor. $\Box$

\begin{corollary}\label{CorNoLeaves}
Let $T$ be a tree such that $\gamma_{\rm m}^\infty(T) =
2\gamma(T)$.  Then no leaf of $T$ belongs to a minimum dominating
set of $T$.
\end{corollary}



\medskip

A neo-colonization $\Pi = \{V_1, V_2, \ldots, V_k\}$ is called
\emph{finest} if it has minimum weight, no parts of size one, and
$k$ is maximum over all such neo-colonizations.

\begin{theorem}\label{ThmIntersect}
Let $D$ be a minimum dominating set of a tree $T$ such that
$\gamma_{\rm m}^\infty = 2\gamma$. If $\Pi = \{V_1, V_2, \ldots,
V_k\}$ is a finest neo-colonization of $T$, then $D \cap V_i \neq
\emptyset$ for $i = 1, 2, \ldots, k$.
\end{theorem}

\noindent{\it Proof.} The proof is by induction on the number of
vertices of $T$.  The statement holds vacuously for all trees with
one vertex.  Suppose it holds for all trees with $n-1$ or fewer
vertices, for some $n \geq 2$.  Let $T$ be a tree with $n$
vertices and $\gamma_{\rm m}^\infty = 2\gamma$.

Let $\Pi = \{V_1, V_2, \ldots, V_k\}$ be a finest neo-colonization
of $T$.  The statement holds if $k = 1$, since $T$ must be a star with at least three vertices.  Hence, assume $k \geq
2$.

Let $x$ be an end vertex of a longest path in $T$.  Then $x$ is a
leaf.  By Corollary \ref{CorNoLeaves}, $x$ is adjacent to a vertex
$y \in D$. Since $\Pi$ is finest, the vertex $x$ and every other
leaf adjacent to $y$ belong to the same part of $\Pi$, say $V_k$.
Let $z$ be the vertex that precedes $y$ on the longest path ending
at $x$.  By Proposition \ref{PropIndep}, the vertex $z$ is not in
$D$.  We distinguish two cases.\\

\noindent {\bf{Case 1.}} Vertex $z$ is a private neighbor of $y$.
Then there is no other stem adjacent to $z$. By choice of $x$ as an
end of a longest path, the vertex $z$ has degree two in $T$,
otherwise there is a longer path.

Suppose $V_k = N[y]$.  Let $T^\prime = T-V_k$.  Then
$\gamma(T^\prime) = \gamma(T) - 1$ and the set $D^\prime = D -
\{y\}$ is a minimum dominating set of $T^\prime$.  Since $V_k$ can
be defended by two guards, $\gamma_{\rm m}^\infty(T^\prime) =
2\gamma(T^\prime)$.  The sequence $\Pi^\prime = \{V_1, V_2, \ldots,
V_{k-1}\}$ is a finest neo-colonization of $T^\prime$.  Hence, by
the induction hypothesis, $D^\prime \cap V_i \neq \emptyset$ for
$i = 1, 2, \ldots, k-1$, and the statement follows.

Suppose $|V_k - N[y]| \geq 2$.  Let $T^\prime = T-N[y]$.  Then
$\gamma(T^\prime) = \gamma(T) - 1$ and the set $D^\prime = D -
\{y\}$ is a minimum dominating set of $T^\prime$.  Since $N[y]$ can
be defended by two guards, $\gamma_{\rm m}^\infty(T^\prime) =
2\gamma(T^\prime)$.  Since both $x$ and $z$ are internal vertices of
the subgraph of $T$ induced by $V_k$, the weight of the
neo-colonization $\Pi^\prime = \{V_1, V_2, \ldots, V_{k-1},
V_k-N[y]\}$ of $T^\prime$ is two less than the weight of $\Pi$. Hence
$\Pi^\prime$ is a finest neo-colonization of $T^\prime$, and the
statement follows from the induction hypothesis as before.


Finally, suppose $V_k = N[y] \cup \{u\}$.    Since $k \geq 2$ and
$z$ is not adjacent to a leaf, the vertex $u$ is adjacent to a
vertex in some other part of $\Pi$, say $V_{k-1}$.  Since the
subgraph of $T$ induced by $V_{k-1}$ is connected, $u$ is adjacent
to at most one vertex in $V_{k-1}$.  Again, let $T^\prime = T-N[y]$.
As before, $\gamma(T^\prime) = \gamma(T) - 1$, the set $D^\prime = D
- \{y\}$ is a minimum dominating set of $T^\prime$, and
$\gamma_{\rm m}^\infty(T^\prime) = 2\gamma(T^\prime)$.   The weight
of the neo-colonization $\Pi^\prime = \{V_1, V_2, \ldots, V_{k-1}
\cup \{u\}\}$ of $T^\prime$ is two less than the weight of $\Pi$.
Hence $\Pi^\prime$ is a finest neo-colonization of $T^\prime$.
Since $z$ is a private neighbor of $y$, we know that the vertex $u
\not\in D$, and the statement follows from the induction hypothesis
as before.\\

\noindent {\bf{Case 2.}} Vertex $z$ is not a private neighbor of
$y$. Then, by Proposition \ref{PropIndep}, $z \not\in D$, and by
Proposition \ref{Prop2EPN}, there is another leaf $w \neq x$
adjacent to $y$.

Suppose first that $V_k = N[y]$ or $V_{k} = N[y] - \{z\}$.   As
before, the tree $T^\prime = T - V_k$ has $\gamma(T^\prime) =
\gamma(T) - 1$, the set $D^\prime = D - \{y\}$ is a minimum
dominating set of $T^\prime$,  $\gamma_{\rm m}^\infty(T^\prime) =
2\gamma(T^\prime)$ and  $\Pi^\prime = \{V_1, V_2, \ldots, V_{k-1}\}$
is a finest neo-colonization of $T^\prime = T-N[y]$.  By the
induction hypothesis, $D^\prime \cap V_i \neq \emptyset$ for $i =
1, 2, \ldots, k-1$, and the statement follows.

Otherwise, $|V_k - (N[y] - \{z\})| \geq 2$.  Then $\Pi^\prime =
\{V_1, V_2, \ldots, V_{k-1}, V_k - (N[y] - \{z\})\}$ is a finest
neo-colonization of $T^\prime = T \setminus (N[y]-\{z\})$.  By the
induction hypothesis, $D^\prime \cap V_i \neq \emptyset$ for $i =
1, 2, \ldots, k-1$, and the statement follows. $\Box$

\begin{corollary}\label{CorPartsBound}
If $\Pi = \{V_1, V_2, \ldots, V_k\}$ is a finest neo-colonization of
a tree $T$ with $\gamma_{\rm m}^\infty = 2\gamma$, then $k \leq
\gamma$.
\end{corollary}

Let $T$ be a tree and $D$ be a minimum dominating set of $T$.  A
\emph{dominating set partition of $T$ with respect to $D$} is a
neo-colonization $\Pi = \{V_1, V_2, \ldots, V_\gamma\}$ such that, for
$i = 1, 2, \ldots, \gamma$, the subgraph of $T$ induced by $V_i$ is
a star centered at a vertex of $D$.  Clearly, every minimum
dominating set $D$ gives rise to at least one dominating set
partition with respect to $D$.  A dominating set partition will be
called \emph{fat} if $|V_i| \geq 3$ for $i = 1, 2, \ldots, \gamma$.

\begin{lemma}\label{LemNeoStars}
Let $T$ be a tree such that $\gamma_{\rm m}^\infty(T) = 2\gamma(T)$.
Then for any minimum dominating set $D$, there exists a fat
dominating set partition of $T$.
\end{lemma}

\noindent{\it Proof.} Let $\Pi = \{V_1, V_2, \ldots, V_\gamma\}$ be
a dominating set partition of $T$ with respect to a minimum
dominating set $D$. Then, by definition, for $i = 1, 2, \ldots,
\gamma$, the subgraph of $T$ induced by $V_i$ is a star centered
at a vertex of $D$. Hence the weight of each part is at most two.
Since the weight of $\Pi$ equals $2\gamma$, and there are $\gamma$
parts, each part has weight exactly two.  Thus each star has at
least three vertices, and $\Pi$ is a fat dominating set partition
of $T$. $\Box$

\begin{lemma} \label{LemFatDP}
Let $T$ be a tree such that $\gamma_{\rm m}^\infty(T) = 2\gamma(T)$.
Then $\Pi = \{V_1, V_2, $ $\ldots, V_\gamma\}$ is a finest
neo-colonization of $T$ in which each part has at least three
vertices if and only if $\Pi$ is a fat dominating set partition with
respect to a minimum dominating set $D$ of $T$.
\end{lemma}

\noindent{\it Proof.} Suppose that $\Pi = \{V_1, V_2, \ldots,
V_\gamma\}$ is a finest neo-colonization of $T$ in which each part
has at least three vertices.  Then each part contributes exactly
two to the weight of $\Pi$.  Thus the subgraph induced by each
part is a tree with connected domination number one, that is, a
star.  If $D$ is the set of center vertices of these stars, then
clearly $D$ is a dominating set of size $\gamma$.

Let $\Pi$ be a  fat dominating set partition with respect to a
minimum dominating set $D$ of $T$.  Then the weight of $\Pi$
equals $2\gamma$, so that $\Pi$ is a minimum weight
neo-colonization of $T$ in which each part has at least three
vertices.  By Corollary \ref{CorPartsBound}, the neo-colonization
$\Pi$ has the maximum number of parts among all minimum weight
neo-colonizations in which there are no parts of size one, hence
$\Pi$ is finest. $\Box$

\begin{theorem}\label{ThmExistChar}
Let $T$ be a tree.  Then $\gamma_{\rm m}^\infty(T) = 2\gamma(T)$ if
and only if there is finest neo-colonization of $T$ which is a fat
dominating set partition with respect to a minimum dominating set of
$T$.
\end{theorem}

\noindent{\it Proof.} The forward implication is immediate from
Lemmas \ref{LemNeoStars} and \ref{LemFatDP}.  For the converse,
suppose there is a finest neo-colonization $\Pi = \{V_1, V_2,
\ldots, V_\gamma\}$ which is a fat dominating set partition with
respect to a minimum dominating set $D$ of $T$.  Then, for $i = 1,
2, \ldots, \gamma$, $|V_i| \geq 3$, the subgraph of $T$ induced by
$V_i$ is a star on at least three vertices.  It follows that the
weight of each part equals two, so that the weight of $\Pi$ equals
$2\gamma$.  Since $\Pi$ has minimum weight, $\gamma_{\rm m}^\infty
= 2\gamma$. $\Box$

\medskip

The definition of a neo-colonization  can be extended to forests.
For each component (tree), the restriction of the neo-colonization
to sets consisting of vertices from that tree, forms a
neo-colonization.

\begin{theorem}\label{ThmLabelCondns}
Let $T$ be a tree.  Then  $\gamma_{\rm m}^\infty = 2\gamma$ if and
only if every minimum dominating set $D$ satisfies
\begin{enumerate}
\item[(a)] No vertex is adjacent to more than two vertices of $D$.
\item[(b)] $D$ is an independent set.
\item[(c)] No two vertices adjacent to two vertices of $D$ are
adjacent.
\item[(d)] Every vertex of $D$ has at least two external private
neighbors.
\item[(e)] There are no two vertices $x, y \in D$ such that the
collection of all private neighbors of $x$ and $y$ induce the
$P_6$ $a, x, b, c, y, z$.
\end{enumerate}
\end{theorem}

\noindent{\it Proof.} Suppose that $\gamma_{\rm m}^\infty =
2\gamma$.  By Proposition \ref{PropIndep}, $D$ is an independent
set.  Hence (b) holds. Suppose there is a vertex $z$ adjacent to at
least three vertices of $D$. Let $\{u, v, w\} \subseteq N(z) \cap
D$.  Form a neo-colonization $\Pi =$ $\{V_1, V_2, \ldots,
V_{\gamma-2}\}$ by letting $V_1 = N[u] \cup N[v] \cup N[w]$ and of
$V_2, V_3, \ldots, V_{\gamma-2}$ be a partition of $V(T) - V_1$ into
stars centered at vertices of $D - \{u, v, w\}$.  The weight of
$V_1$ is five, so the weight of $\Pi$ is at most $5+2(\gamma - 3) <
2\gamma$. Since $\gamma_{\rm m}^\infty = 2\gamma$, it follows that
(a) holds.

Suppose that $x$ and $y$ are adjacent vertices that are each
adjacent to two vertices of $D$.   Let $N[x] \cap D = \{a, b\}$ and
$N[y] \cap D = \{c, d\}$.  Then $\{a, b\} \cap \{c, d\} =
\emptyset$. Form a neo-colonization $\Pi = \{V_1, V_2, \ldots,
V_{\gamma-3}\}$ by letting $V_1 = N[a] \cup N[b] \cup N[c] \cup N[d]$
and $V_2, V_3, \ldots, V_{\gamma-3}$ be a partition of $V(T) - V_1$
into stars centered at vertices of $D - \{a, b, c, d\}$. The weight
of $V_1$ is seven, so that the weight of $\Pi$ is at most
$7+2(\gamma - 4) < 2\gamma$. Since $\gamma_{\rm m}^\infty =
2\gamma$, it follows that condition (c) holds.

Suppose there is a vertex $x \in D$ with only one external private
neighbor, say $y$,
 Form a neo-colonization  $\Pi = \{V_1, V_2, \ldots, V_{\gamma-1}, \{x, y\}\}$ such that, for $i = 1, 2, \ldots, \gamma-1$, the subgraph of $T$ induced by $V_i$ is a star centered at a vertex of $D-\{x\}$. Then the weight of $\Pi$ is at most $2(\gamma - 1) + 1 < 2\gamma$. Since $\gamma_{\rm m}^\infty = 2\gamma$, it follows that the component $V_1$ can not exist.  Therefore, condition (d) holds.

Suppose that $x, y \in D$ and the set of all private neighbors of
$x$ and $y$ comprise the path $x_1, x, x_2, y_1, y, y_2$.  Form a
neo-colonization $\Pi = \{V_1, V_2,$ $ \ldots, V_{\gamma-2},$ $\{x_1,
x\}, \{x_2, y_1\}, \{y, y_2\}\}$, where $V_1, V_2, \ldots,
V_{\gamma-2}$ is a partition of $V(T) - $ 
$\{x, x_1, x_2, y, y_1, y_2\}$ into stars centered at vertices of
$D - \{x, y\}$. Then the weight of $\Pi$ is at most $2(\gamma - 2)
+ 3 < 2\gamma$. Since $\gamma_{\rm m}^\infty = 2\gamma$, it
follows that the vertices $x$ and $y$ can not exist.  Therefore,
condition (e) holds.

The proof of the converse implication is by induction on the
number of vertices of $T$.  The statement holds for all trees on
one vertex.  Suppose it holds for all trees on at most $n-1$
vertices, for some $n \geq 2$.  Let $T$ be a tree on $n$ vertices
for which conditions (a) through (e) hold.  We need only consider
the case when $T$ is not a star.  Let $D$ be a minimum dominating set of $T$.

Let $x$ be an end vertex of a longest path of $T$, and $y$ the
unique neighbor of $x$.   By choice of $x$, the vertex $y$ has a
unique neighbor $z$ which is not a leaf.
Since condition (d) holds,
no leaf is in $D$. Hence $y\in D$. Then, by (b),
the vertex $z$ can not belong to any minimum dominating set.
We claim that  the vertex $z$ has
degree at most three.  Suppose $z$ has four neighbors, $y, u, v, w$, where $y$ and $w$ lie on a longest path ending at $x$.  By choice of $x$, the vertices $u$ and $v$ are either leaves or adjacent to a leaf.  By condition (d), no leaf is in $D$.  Hence each of $u$ and $v$ is adjacent to a leaf, and $u, v \in D$. Since $y \in D$, the vertex $z$ has three neighbors in $D$, contrary to (a).  Thus, $2 \leq \mathit{deg}_T(z) \leq 3$.  We
consider these cases separately.\\

\noindent {\bf{Case 1.}} $\mathit{deg}_T(z) = 3$. Let $N_T(z) = \{a,
b, y\}$, where $b$ lies on a longest path starting at $x$.  Then $a$ is not a leaf,
but every vertex in
$N_T(a) - \{z\}$ is a leaf.  Since condition (d) holds, the vertices
$a$ and $y$ are each adjacent to at least two leaves.

Let $T^\prime = T - (N[y] - \{z\})$.  Then $\gamma(T^\prime) =
\gamma(T) - 1$, and $D^\prime$ is a minimum dominating set of
$T^\prime$ if and only if $D^\prime \cup \{y\}$ is a minimum
dominating set of $T$.  Since condition (a) holds, no minimum
dominating set of $T$ contains $b$.  Therefore there is no minimum
dominating set of $T^\prime$ for which the vertex $z$ belongs to a
$P_6$ as described in condition (e). It follows that conditions
(a) through (e) therefore hold for every minimum dominating set of
$T^\prime$. By the induction hypothesis, $\gamma_{\rm
m}^\infty(T^\prime) = 2 \gamma(T^\prime)$.

Let $\Pi = \{V_1, V_2, \ldots, V_k\}$ be a minimum weight
neo-colonization of $T$.   All leaves adjacent to $y$ belong to the
same part as $y$, say $V_k$ (otherwise $\Pi$ does not have minimum
weight).  Since the weight of $\Pi$ is at most $2\gamma(T)$, it
suffices to show that it is also at least $2\gamma(T)$.  There are
three possibilities, depending on $V_k - (N[y] - \{z\})$.

If $V_k = N[y] - \{z\}$, then $\Pi^0 = \{V_1, V_2, \ldots,
V_{k-1}\}$ is a neo-colonization of $T^\prime$, and hence has
weight at least $2(\gamma(T) - 1)$.  Therefore the weight of $\Pi$
is at least $2\gamma(T)$.

If $V_k = N[y]$, then the weight of $\Pi$ is one more than the
weight of the neo-colonization $\Pi^1 = \{V_1, V_2, \ldots,
V_{k-1}, \{z\}\}$ of $T^\prime$.  Since all leaves adjacent to $a$
must belong to the same part as (a), and since the subgraph
induced by each part is connected, $N[a] - \{z\}$ is a part of
$\Pi^1$, say $V_{k-1} = N[a] - \{z\}$.  The neo-colonization
$\{V_1, V_2, \ldots, V_{k-1} \cup \{z\}\}$ has weight one less that
$\Pi^1$, so that the weight of $\Pi^1$ is at least $2(\gamma(T) -
1) + 1$.  Therefore the weight of $\Pi$ is at least $2\gamma(T)$.

If $V_k - N[y] \neq \emptyset$, then the connected domination number
of the subgraph of $T$ induced by $V_k$ is two more than the
connected domination number of the subgraph induced by $V_k-(N[y] -
\{z\})$.  Hence, the weight of $\Pi$ is two more than the
neo-colonization $\Pi^2 = \{V_1, V_2, \ldots, V_{k-1}, V_k-(N[y] -
\{z\})\}$ of $T^\prime$ of $T^\prime$.  Since the weight of $\Pi^2$
is at least $2(\gamma(T) - 1)$, the weight of $\Pi$ is at least
$2\gamma(T) $.\\

\noindent {\bf{Case 2.}} $\mathit{deg}_T(z) = 2$. Let $N_T(z) = \{b,
y\}$, and $T^\prime = T - N[y]$.  Then $\gamma(T^\prime) = \gamma(T)
- 1$, and $D^\prime$ is a minimum dominating set of $T^\prime$ if
and only if $D^\prime \cup \{y\}$ is a minimum dominating set of
$T$.  Since conditions (a) through (e) hold for every minimum
dominating set of $T$, they also hold for every minimum dominating
set of $T^\prime$. Therefore, $\gamma_{\rm m}^\infty(T^\prime) = 2
\gamma(T^\prime)$.

Let $D^\prime$ be a minimum dominating set of $T^\prime$. Since
condition (d) holds for $T$, the vertex $b \not\in D^\prime$.

Let $\Pi = \{V_1, V_2, \ldots, V_k\}$ be a minimum weight
neo-colonization of $T$. Then all leaves adjacent to $y$ belong to
the same part as $y$, say $V_k$ (otherwise $\Pi$ does not have
minimum weight).  If $z$ also belongs to $V_k$, the weight of
$\Pi$ is two more than the weight of $\Pi^\prime$, and it follows
that $\gamma_{\rm m}^\infty(T) = 2 \gamma(T)$.

Hence, suppose that $z$ does not belong to $V_k$.  If there are at
least two leaves adjacent to $y$, or if the part  containing $z$,
say $V_{k-1}$, is not a star centered at $b$, then the weight of
$\Pi$ is two more than the weight of $\Pi^\prime$, and it follows
that $\gamma_{\rm m}^\infty(T) = 2 \gamma(T)$.

Assume then, that $x$ is the only leaf adjacent to $y$ and the part
$V_{k-1}$ containing $z$ is a star centered at $b$.   The list
$\Pi^{\prime\prime} = \{V_1, V_2, \ldots, V_{k-2}\}$ is a
neo-colonization of $T^\prime - V_{k-1}$ which is of weight at most
$2(\gamma(T)-1) - 1$.  Hence the domination number of $T^\prime -
V_{k-1}$ is at most $\gamma(T) - 2 < \gamma(T^\prime)$.  Since
$V_{k-1}$ is a star centered at $b$, there is a dominating set
$D^{\prime\prime}$ of $T^\prime$ that contains $b$.  Therefore
$D^{\prime\prime} \cup \{y\}$ is a dominating set of $T$ in which
$y$ has only one external private neighbor.  Hence condition (d)
does not hold for $T$, a contradiction.

The result now follows by induction. $\Box$

\medskip

Let $D$ be a dominating set of the tree $T$.  The \emph{domination
labeling of $T$ with respect to $D$} is the function $\ell_D: V \to
\{1, 2, \ldots, |D|\}$ that assigns to each vertex $x \in V$ the
integer $\ell_D(x) = |N[x] \cap D|$.

Given a domination labeling $\ell_D$, we use $F_1(T)$, or $F_1$ when
the context is clear, to denote the subgraph of $T$ induced by the
set of vertices that are labeled one.

\begin{corollary}
Let $T$ be a tree.  Then  $\gamma_{\rm m}^\infty = 2\gamma$ if and
only if every minimum dominating set $D$ satisfies
\begin{enumerate}
\item[(a)] $1 \leq \ell_D(v) \leq 2$ for every $v \in V$.
\item[(b)] $\ell_D(x) = 1$ for every $x \in D$. \item[(c)] The set
$L_2 = \{x: \ell_D(x) = 2\}$ is an independent set. \item[(d)] If
$x \in D$, then $|N_{F_1}(x)| \geq 2$.
\item[(e)] There are no two vertices $x, y \in D$ such that the
subgraph of $F_1$ induced by $N_{F_1}[x] \cup N_{F_1}[y]$ is
isomorphic to $P_6$.
\end{enumerate}
\end{corollary}

\begin{corollary}
Let $T$ be a tree.  Then  $\gamma_{\rm m}^\infty < 2\gamma$ if and
only if there exists a minimum dominating set $D$ and a domination
partition with respect to $D$ such that at least one of the
following statements holds:
\begin{enumerate}
\item[(a)] Some vertex is adjacent to three vertices of $D$.
\item[(b)] Two vertices of $D$ are adjacent. \item[(c)] There
are adjacent vertices $x$ and $y$ that are each adjacent to two
vertices of $D$. \item[(d)] The subgraph induced by some part of the
partition is isomorphic to $K_1$ or $K_2$. \item[(e)] There are two
parts $X$ and $Y$ such that the subgraph of $T$ induced by $X \cup
Y$ is isomorphic to $P_6$.
\end{enumerate}
\end{corollary}

\section{Trees with $\gamma_{\mathrm{m}}^\infty = \beta$}

The following is a fundamental upper bound on the $m$-eternal
domination number.

\begin{theorem} \label{beta-bound}
\emph{\cite{GHH}}\label{alpha-bound} Let $G$ be a graph. Then
$\gamma_{\mathrm{m}}^{\infty}(G) \leq \beta(G)$.
\end{theorem}

Let $s$ be a stem in a tree $T$ with at least two vertices.  We call $s$ \emph{exposed} if it has at most one neighbor that is an internal vertex of $T$. Note that a tree with at most two vertices does not have exposed stems by definition.

The next two propositions have simple proofs which we omit.

\begin{proposition}
Every tree with at least two vertices that is not a star has at least two exposed stems.
\label{TwoStrong}
\end{proposition}

\begin{proposition}
Iteratively deleting exposed stems and the leaves adjacent to them partitions the vertices of a tree $T$ into subsets that each induce a star, and $\beta(T)$ is the sum of the independence numbers of the stars in the partition.
\label{BetaTree}
\end{proposition}

We consider the following operation which is a restriction of the procedure in Proposition \ref{BetaTree}.\\

\noindent {\bf Operation EWS:} \emph{If there is an exposed weak stem, delete it and its leaf.}

\begin{lemma}{\rm \cite{KM2}}
If $T^\prime$ results from one application of operation EWS to $T$,  then $\beta(T^\prime) = \beta(T)-1$ and $\gamma_{\rm m}^\infty(T^\prime) = \gamma_{\rm m}^\infty(T)-1$.
\label{R1}
\end{lemma}

The following lemma can be proved trivially by induction.

\begin{lemma}
If $T$ has an exposed strong stem that is adjacent to more than two leaves, then $\gamma_{\rm m}^\infty(T) < \beta(T)$.
\label{TooManyLeaves}
\end{lemma}

%

\begin{theorem} A tree with at least two vertices has $\beta = \gamma_{\rm m}^\infty$ if and only if repeated applications of operation EWS reduces $T$ to $K_1$ or $K_2$.
\end{theorem}

\noindent{\it Proof.}
Suppose $k$ applications of operation EWS reduces $T$ to $K_1$ or $K_2$.  Form a neo-colonization $\{V_1, V_2, \ldots, V_{k+1}\}$ of $T$, where $V_i$ is the set containing the stem and leaf deleted by the $i^{th}$ application of the operation, $1 \leq i \leq k$, and $V_{k+1}$ is the vertex set of the final $K_1$ or $K_2$.
By the discussion above, the weight of this neo-colonization is $\beta$.  We show that there is no neo-colonization of smaller weight by induction on $n = |V|$.

By inspection, the statement is true for trees with two or three vertices.  Suppose it holds for trees with between two and $n-1$ vertices, for some $n \geq 4$.  Let $T$ be a tree with $n$ vertices.  If operation EWS can be applied to $T$, then the result follows by induction.  Hence suppose EWS cannot be applied to $T$.  Then every exposed stem is strong.  Further, by Lemma \ref{TwoStrong}, $T$ has at least two strong stems.

By Lemma \ref{TooManyLeaves}, if there is an exposed strong stem which is adjacent to more that two leaves then $\beta(T) > \gamma_{\rm m}^\infty$ and the statement follows.  Hence assume every exposed strong stem is adjacent to exactly two leaves.  Since $T$ has at least four vertices, it follows that every exposed strong stem is adjacent to a unique internal vertex of $T$.

Let $s$ be an exposed strong stem of $T$ with maximum eccentricity.  Let $X$ be the set consisting of $s$ and the two leaves to which it is adjacent.  It is clear that $\beta(T) = 2 + \beta(T-X)$.  We now consider the outcome of applying operation EWS to $T-X$.

Suppose first that operation EWS does not reduce $T-X$ to $K_1$ or $K_2$.  Then, by the induction hypothesis, $\beta(T-X) > \gamma_{\rm m}^\infty(T-X)$.   Let $\{V_1, V_2, \ldots, V_p\}$ be a minimum weight neo-colonization of $T-X$. Then $\{V_1, V_2, \ldots, V_p, X\}$ is a neo-colonization of $T$ of weight $2 + \gamma_{\rm m}^\infty(T-X)$.  In this case we have
\begin{eqnarray*}
\beta(T) & = & 2  + \beta(T-X)\\
& > & 2 + \gamma_{\rm m}^\infty(T-X)\\
& \geq & \gamma_{\rm m}^\infty(T),
\end{eqnarray*}
and the result follows.

Now suppose that operation EWS reduces $T-X$ to $K_1$ or $K_2$.  We claim that it is reduced to $K_1$.  Since $T$ has at least two exposed strong stems, $T-X$ has at least one, say $w$.  The reduced graph arising from each application of EWS is a tree.  No application of EWS deletes $w$ or a leaf adjacent to $w$ unless the tree under consideration is a path on three vertices.  This proves the claim.

Next, we claim that $T$ is a caterpillar with spine $s, x_1, x_2, ..., x_k, w$, such that
\begin{enumerate}
\item[(a)] each of $x_1, x_2, ..., x_k$ is either a loner or a weak stem.
\item[(b)] there are two leaves adjacent to $s$ and to $w$;
\end{enumerate}
To see that $T$ is a caterpillar, first recall from above that EWS cannot be applied to $T$. Further observe that $T$ can have no vertex $v$ adjacent to three internal vertices because EWS would delete all internal vertices on the path from $s$ to $v$ in $T-X$ and then could not be applied again. (Since $w$ is an exposed stem, $v \neq w$.) Similarly, neither $T-X$ nor any tree derived from it by applying operation EWS has at least two strong stems.  Point (a) now follows.   Point (b) was established above.  This proves the claim.

We complete the proof by showing that $\beta(T) > \gamma_{\rm m}^\infty(T)$.  Partition the spine of $T$, i.e., the vertices  $s, x_1, x_2, ..., x_k, w$, into maximal paths of stems and maximal paths of loners.  Since none of the reduced trees arising from $T-X$ has two strong stems, and $T-X$ is reduced to $K_1$, each maximal path of loners has an even number of vertices.


Suppose first that $T$ has no loners.  Then $T$ has connected domination number $k+2$ and independence number $k+4$.  Hence assume $T$ has loners.  Construct a neo-colonization of $T$ by forming a part out of each maximal path of consecutive stems along the spine of $T$, the leaves adjacent to them, and any loner adjacent to exactly one of these stems.  The remaining vertices are all loners belonging to disjoint paths on an even number of vertices.  Partition each of these paths into $K_2$'s.

We argue that there is an independent $I$ set of size greater than the weight of the neo-colonization.  The part containing $s$ has independence number two more than connected domination number.  Put these vertices into $I$.  The maximum independent set in this part contains the vertex that was a loner in $T$.  Each $K_2$ part has weight one and contributes one vertex to $I$.  The vertex that it contributes is the one farthest from $s$ along the spine of $T$.  Proceeding away from $s$ along the spine of $T$, eventually there is a part containing a maximal path of stems.  It has independence number two more than connected domination number, but the loner closest to $s$ in this part is cannot be included in $I$ because it is adjacent to the loner from the previous part which belongs to $I$.  Put the remainder of the maximum independent set of the part, i.e. all but this loner,  into $I$.  Continuing in this way, eventually the part containing $w$ contributes its connected domination number plus one vertices to $I$.  Thus the weight of the neo-colonization is $|I| - 1$. This completes the proof.
$\Box$

\medskip

An alternate characterization is given as a corollary of the previous result.

\begin{corollary}\label{ind-char2}
Let $T$ be a tree with at least two vertices. Then $\gamma_{\mathrm{m}}^{\infty}(T) =
\beta(T)$ if and only if there exists a minimum-weight
neo-colonization of $T$ containing only parts that are $K_2$'s or
$P_3$'s and there is at most one $P_3$.
\end{corollary}

$K_{1, 3}$ with each edge subdivided twice is an example showing
that the minimum weight condition on neo-colonization is
necessary in Proposition \ref{ind-char2}. This graph has a
neo-colonization consisting of $K_2$'s and $P_3$'s but has
$\gamma_{\mathrm{m}}^{\infty }(T) < \beta(T)$.

\section{Concluding Remarks}

We begin this section by stating a result of Chambers et al. \cite{prince2} which also appears in the survey \cite{KM8}.

\begin{theorem}\cite{prince2, KM8}
For all connected graphs, $\gamma_{\mathrm{m}}^\infty \leq \lceil \frac{|V|}{2} \rceil$.
\end{theorem}

The trees for which $\gamma_{\mathrm{m}}^\infty = \gamma$ and those for which $\gamma_{\mathrm{m}}^\infty = \beta$ achieve equality in this bound.  There are other trees for which equality holds.  One example is the tree obtained from the path $v_1, v_2, \ldots, v_9$ by adding two vertices $v_{10}, v_{11}$ and the edges $v_4v_{10}, v_6v_{11}$.  Characterizing the trees such that $\gamma_{\mathrm{m}}^\infty = \lceil \frac{|V|}{2} \rceil$ remains an open problem.

\bigskip\bigskip
\noindent{\large\bf{Acknowledgments}}\\
We thank an anonymous referee for their careful reading of the
paper, their helpful comments and corrections.


\begin{thebibliography}{99}                                                                                               %


\bibitem {ABB}M.~Anderson, C.~Barrientos, R.~Brigham, J.~Carrington,
R.~Vitray, J.~Yellen, Maximum demand graphs for eternal security, $\emph{J}%
$.~\emph{Combin.~Math.~Combin.~Comput.}~\textbf{61} (2007), 111--128.


\bibitem {BCG2}A.P.~Burger, E.J.~Cockayne, W.R.~Gr\"{u}ndlingh, C.M.~Mynhardt,
J.H.~van Vuuren, W.~Winterbach, Infinite order domination in graphs,
\emph{J.~Combin.~Math.~Combin.~Comput.}~\textbf{50} (2004), 179--194.

\bibitem{Fink} J.F.~Fink, M.S.~Jacobson, L.F.~Kinch and J.~Roberts,
On graphs having domination number half their order.
\emph{Period. Math. Hungar.}~\textbf{16} (1985), 287--293.

\bibitem{Payan}
C.~Payan and N.H.~Xuong,
Domination-balanced graphs.
\emph{J. Graph Theory} \textbf{6} (1982), 23--32.

\bibitem{prince2} E. Chambers, W. Kinnersly, and N. Prince, Mobile eternal security in graphs, manuscript (2008).






\bibitem {GHH}W.~Goddard, S.M.~Hedetniemi, S.T.~Hedetniemi, Eternal security
in graphs, \emph{J.~Combin.~Math.~Combin.~Comput.}~\textbf{52} (2005), 169--180.

\bibitem {GK}J.~Goldwasser, W.F.~Klostermeyer, Tight bounds for eternal
dominating sets in graphs, \emph{Discrete Math.}~\textbf{308} (2008), 2589--2593.


\bibitem {HHS}T.~W.~Haynes, S.~T.~Hedetniemi, P.~J.~Slater, \emph{Fundamentals
of Domination in Graphs}. Marcel Dekker, New York, 1998.





\bibitem {KM}W.F.~Klostermeyer, G.~MacGillivray, Eternal security in graphs of
fixed independence number, $\emph{J}$.~\emph{Combin.~Math.~Combin.~Comput.}%
~\textbf{63} (2007), 97--101.

\bibitem {KM2}W.F.~Klostermeyer, G.~MacGillivray, Eternal dominating sets in
graphs, $\emph{J}$.~\emph{Combin.~Math.~Combin.~Comput.}~\textbf{68} (2009), 97--111.


\bibitem {KM6}W.F.~Klostermeyer and C.M.~Mynhardt, Graphs with Equal Eternal Vertex Cover and
Eternal Domination Numbers, {\it Discrete Math.} ~\textbf{311}
(2011), 1371--1379.

\bibitem {KM7}W.F.~Klostermeyer and C.M.~Mynhardt, Vertex Covers
and Eternal Dominating Sets, {\it Discrete Applied Mathematics} ~\textbf{160} (2012), pp. 1183-1190.

\bibitem {KM8}W.~F.~Klostermeyer, C.~M.~Mynhardt,  Protecting a Graph with Mobile Guards, to appear in {\bf{ Movement on networks}}, Cambridge University Press.
\end{thebibliography}
\end{document}